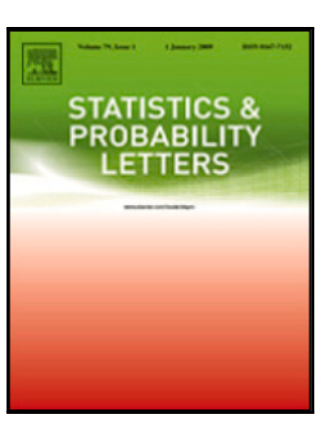

# A new and flexible class of sharp asymptotic time-uniform confidence sequences

Felix Gnettner *, Claudia Kirch

*Otto-von-Guericke-Universität Magdeburg, Fakultät für Mathematik,*
*Institut für Mathematische Stochastik, Universitätsplatz 2, 39106 Magdeburg, Germany*

ABSTRACT

Confidence sequences are anytime-valid analogues of classical confidence intervals that do not suffer from multiplicity issues under optional continuation of the data collection. As in classical statistics, asymptotic confidence sequences are a nonparametric tool showing under which high-level assumptions asymptotic coverage is achieved so that they also give a certain robustness guarantee against distributional deviations. In this paper, we propose a new flexible class of confidence sequences yielding sharp asymptotic time-uniform confidence sequences under mild assumptions. Furthermore, we highlight the connection to corresponding sequential testing problems and detail the underlying limit theorem.

## 1. Introduction

While a classical $(1-\alpha)$-confidence interval $\widetilde{C}_t = \widetilde{C}(X_1, \dots, X_t)$ for the expected value $\mu_X$ of iid observations $X_1, \dots, X_t$ satisfies

$$\mathbb{P}(\mu \in \widetilde{C}_t) \geqslant 1 - \alpha \quad \text{for all } t \in \mathbb{N},$$

a time-uniform confidence sequence $C_t = C(X_1, \dots, X_t), t \in \mathbb{N}$, for $\mu_X$ fulfils the stronger property

$$\mathbb{P}(\mu \in C_t \text{ for all } t \in \mathbb{N}) \geqslant 1 - \alpha.$$

The latter property allows repeated calculation and updating of the confidence bounds as new data arise without violating the statistical coverage guarantees. Such sequences have already been proposed as early as the 1970s (Robbins, 1970; Lai, 1976), but have recently attracted broad attention in the context of safe testing and inference methods (Grünwald et al., 2024). While most of the earlier proposals for confidence sequences are based on a parametric model, Robbins and Siegmund (1970) already provided an asymptotic nonparametric construction technique for iid centred random variables with known variance.

Recently, Waudby-Smith et al. (2024) proposed a general construction principle for deriving asymptotic confidence sequences based on strong invariance principles. Bibaut et al. (2024) consider corresponding sequential tests and provide asymptotic type-I-error as well as expected-rejection-time guarantees under general nonparametric data generating processes. Moreover, Bibaut et al. (2021) construct an asymptotic nonparametric confidence sequence by means of a weak invariance principle for martingale difference triangular arrays. A detailed overview of recent developments is provided by Ramdas et al. (2023).

In this paper, we propose a new flexible construction based on different asymptotic tools that can be achieved under somewhat weaker assumptions and have proven to be useful in the context of monitoring data for change points (see the discussion in

* Corresponding author.
*E-mail address:* felix.gnettner@ovgu.de (F. Gnettner).

Section 5). First, we focus on uniform confidence sequences for the location problem in Section 2, then highlight the connection to sequential testing in Section 3, before we detail the limit theorem behind these constructions in Section 4.

## 2. Sharp uniform confidence sequences for a location problem

Waudby-Smith et al. (2024, Definition 2.7) define sharp asymptotic $(1-\alpha)$-confidence sequences $(C_t(m;\alpha))_{t\geqslant m}$ for a parameter $\mu$ as sequences satisfying $\lim_{m\to\infty}\mathbb{P}\left(\mu\in C_t(m;\alpha)\quad\text{for all } t\geqslant m\right)=1-\alpha$. The parameter $m$ plays the same role here as the sample size in classical statistics and can be thought of as a burn-in period for the asymptotic approximation to work well. As an example, Waudby-Smith et al. (2024, Theorem 2.8) prove that $C_t(m;\alpha)=\hat{\mu}_t\pm\hat{\sigma}_t\,\tilde{b}_{t,m}(\alpha)$ is a sharp asymptotic confidence sequence for the mean of iid random variables $(X_i)_{i\in\mathbb{N}}$ with finite variance $\sigma_X^2>0$, where $\hat{\mu}_t$, $\hat{\sigma}_t^2$ are the sample mean and sample variance based on the first $t$ observations. Their $\tilde{b}_{t,m}(\alpha)$ has a specific shape inherited from a boundary crossing result for Wiener processes. Indeed, the shape of the boundary sequence $\tilde{b}_{t,m}(\alpha)$ determines the length of the confidence interval at the time point $t$. It is connected with how quickly the $\alpha$ in the confidence level is *spent* across the time points $t$ (with the words of Lan and DeMets (1983)) while keeping uniform coverage across the sequential data collection. For some applications, it might be more desirable to start with smaller confidence intervals even if they are larger compared to different choices for the boundary curves at later times and vice versa for others. As such, some flexibility in choosing the boundary curves $\tilde{b}_{t,m}(\alpha)$ is desirable.

The following construction is based on a different type of limit result (see Section 4) allowing for a lot of flexibility in the choice of the boundary curves $\tilde{b}_{t,m}(\alpha)$. Additionally, this underlying limit theorem combines a functional central limit result with finite sample concentration inequalities, permitting the construction of confidence sequences that hold the level uniformly for all $t\in\mathbb{N}$ and not just for all $t\geqslant m$ as in Waudby-Smith et al. (2024, Definition 2.7):

**Definition 2.1.** A sequence of intervals $(C_t(m;\alpha))_{t\in\mathbb{N}}$ based on $(X_t)_{t\in\mathbb{N}}$ is a sharp asymptotic $(1-\alpha)$-confidence sequence for $\mu$ if it satisfies $\lim_{m\to\infty}\mathbb{P}\left(\mu\in C_t(m;\alpha)\quad\text{for all } t\in\mathbb{N}\right)=1-\alpha$ for any $\alpha\in(0,1)$.

Previously proposed asymptotic confidence sequences are included in the above definition by using the trivial choice $C_t(m)=\mathbb{R}$ for all $t<m$. Throughout the paper, we make the following assumptions on the statistical model:

**Assumption $\mathcal{M}$.** Let $(X_t)_{t\in\mathbb{N}}$ be a sequence of iid random variables with mean $\mu_X$ and variance $0<\sigma_X^2<\infty$. Let $\hat{\sigma}_t^2\to\sigma_X^2$ a.s. with $\hat{\sigma}_t^2>0$ a.s. for all $t\in\mathbb{N}$.

We use the assumption of independent errors for ease of presentation only. Indeed, all results hold analogously for stationary time series under weak nonparametric assumptions by replacing the variance with the long-run variance, see Remark 4.2. In practice, confidence sequences should only be given when the chosen variance estimator is strictly positive. While for continuous random variables this is fulfilled as soon as $t\geqslant 2$, this is not necessarily true for discrete random variables, so that additional adaptations are necessary. Formally, this can be dealt with in the above framework by replacing the variance estimator with $\infty$ until the original estimator is strictly positive, as such an adaptation does not change the strong consistency of the sequence.

We propose to use the following sequence of open confidence intervals

$$C_t(m;\alpha)=\hat{\mu}_t\pm\hat{\sigma}_t\cdot c_\alpha(\rho)\cdot b_t(m;\rho),\qquad\text{where}\quad \hat{\mu}_t=\frac{1}{t}\sum_{i=1}^{t}X_i,\quad b_t(m;\rho)=\frac{\sqrt{m}}{t\,\rho(t/m)}. \tag{1}$$

The function $\rho$ determines the shape of the boundary curve and can be chosen flexible as long as it satisfies the following Assumption $\mathcal{BC}$, while $c_\alpha(\rho)$ is a constant depending on $\alpha$ and $\rho$ chosen such that the sequence is sharp; see Theorem 2.2 below.

**Assumption $\mathcal{BC}$.** Let $\rho:(0,\infty)\to[0,\infty)$ be a mapping that is continuous if restricted to $(0,e_\rho)$, where $e_\rho=\sup\{s>0:\rho(s)>0\}$. Additionally, the following conditions hold

$$(A_1)\ \limsup_{s\to 0}s^{\gamma_1}\cdot\rho(s)<\infty\ \text{for some}\ \gamma_1\in[0,1/2),\qquad (A_2)\ \limsup_{s\to\infty}s^{1-\gamma_2}\cdot\rho(s)<\infty\ \text{for some}\ \gamma_2\in[0,1/2).$$

Typically, the endpoint $e_\rho$ will be equal to infinity, corresponding to an open-end procedure where data are collected possibly forever, while for a finite endpoint at most $\lfloor m\,e_\rho\rfloor$ data points are collected.

The following theorem is a direct consequence of Theorem 4.1 in Section 4.

**Theorem 2.2.** *Under Assumptions $\mathcal{M}$ and $\mathcal{BC}$ the sequence $C_t(m;\alpha)$ as in (1) is a sharp asymptotic $(1-\alpha)$-confidence sequence in the sense of Definition 2.1, if $c_\alpha(\rho)$ is chosen as the $(1-\alpha)$-quantile of $\sup_{y>0}|\rho(y)\cdot W(y)|$, where $(W(y))_{y\in(0,\infty)}$ denotes a standard Wiener process.*

In practice, we propose the use of the weight functions $\rho(s)=(1+s)^{\gamma_1+\gamma_2-1}/s^{\gamma_1}$, for some $0\leqslant\gamma_1,\gamma_2<1/2$, because for this class the value $c_\alpha(\rho)$ can be chosen as the $(1-\alpha)$-quantile of $\sup_{0\leqslant y\leqslant 1}\frac{|B(y)|}{y^{\gamma_1}(1-y)^{\gamma_2}}$, where $(B(y))_{y\in[0,1]}$ denotes a standard Brownian bridge (see Proposition 4.3). In contrast to a general weight function, the supremum is now taken over a finite interval and the necessary quantiles can e.g. be approximated as in Franke et al. (2022). For a choice of $\gamma_2>0$ the length of the confidence sequences converges to zero (for every fixed $m$ and $t\to\infty$), an observation connected to Theorem 3.2.

## 3. Connection to sequential testing problems and safe testing

Analogously to the well-known duality between classical statistical tests and confidence intervals, there is a close connection between confidence sequences and sequential testing. Indeed, the above approach is closely related to methodology well-established in nonparametric sequential change point testing first proposed by Chu et al. (1996), Horváth et al. (2004); see Aue and Kirch (2024) for a recent survey on the subject.

To elaborate, consider testing $H_0(\mu)\ :\ \mu_X = \mu$ against $H_1(\mu)\ :\ \mu_X \neq \mu$ in the situation of the previous section. The above confidence sequence is obtained by inverting the sequence of tests $(\varphi_t(m;\mu,\alpha))_{t\in\mathbb{N}}$ that reject at time point $t$ if $|\widehat{\mu}_t - \mu| \geqslant \widehat{\sigma}_t \cdot c_\alpha(\rho)\cdot b_t(m;\rho)$ and vice versa. Therefore, sequentially testing in this manner can be used in the context of safe testing as it controls the family-wise error rate (FWER) across $t$ at level $\alpha$.

The following theorem is equivalent to Theorem 2.2 and a direct consequence of Theorem 4.1 in Section 4.

**Theorem 3.1.** *Consider the sequence of tests* $\varphi_t(m;\mu,\alpha) = \mathbb{1}\left\{|\widehat{\mu}_t - \mu| \geqslant \widehat{\sigma}_t \cdot c_\alpha(\rho)\cdot b_t(m;\rho)\right\}$ *with the notation and under the assumptions of Theorem 2.2. Then, for any* $\alpha \in (0,1)$ *and any* $\mu \in \mathbb{R}$, *it holds*

$$\lim_{m\to\infty} \mathbb{P}_{\mu_X=\mu}\left(\varphi_t(m;\mu,\alpha) = 1 \text{ for some } t \in \mathbb{N}\right) = \alpha.$$

Furthermore, the above sequence of tests has the following finite-sample stopping guarantee.

**Theorem 3.2.** *Under the assumptions of Theorem 3.1 consider* $\rho(\cdot)$ *with* $\liminf_{s\to\infty} s^{1-\tilde{\gamma}_2}\cdot\rho(s) > 0$ *for some* $\tilde{\gamma}_2 \in (0, 1/2)$*, then for any* $m$ *and any* $\mu_X \neq \mu$ *it holds*

$$\mathbb{P}_{\mu_X\neq\mu}\left(\varphi_t(m;\mu,\alpha) = 1 \text{ for some } t\in\mathbb{N}\right) = 1.$$

*If* $\rho(\cdot)$ *is chosen such that the above condition only holds with* $\tilde{\gamma}_2 = 0$*, then:* $\lim_{m\to\infty}\mathbb{P}_{\mu_X\neq\mu}\left(\varphi_t(m;\mu,\alpha) = 1 \text{ for some } t\in\mathbb{N}\right) = 1.$

**Proof.** By the law of large numbers $\lim_{t\to\infty}|\widehat{\mu}_t - \mu| > 0$ a.s., while $\widehat{\sigma}_t\cdot c_\alpha(\rho)\cdot b_t(m;\rho)\to 0$ a.s. by assumption as

$$b_t(m;\rho) = \frac{\sqrt{m}}{t\,\rho(t/m)} = \frac{1}{t^{\tilde{\gamma}_2}\, m^{1/2-\tilde{\gamma}_2}}\;\frac{1}{(t/m)^{1-\tilde{\gamma}_2}\,\rho(t/m)}. \qquad \square$$

Furthermore, it is possible to control the family-wise error rate (FWER) when simultaneously testing one-sided hypotheses for $\mu = \mu_1,\ldots,\mu_s$ of the form

$$H_0^{(r)}(\mu)\ :\ \mu_X \leqslant \mu \quad\text{against}\quad H_1^{(r)}(\mu)\ :\ \mu_X > \mu, \qquad H_0^{(\ell)}(\mu)\ :\ \mu_X \geqslant \mu \quad\text{against}\quad H_1^{(\ell)}(\mu)\ :\ \mu_X < \mu \tag{2}$$

with corresponding sequential one-sided tests

$$\varphi_t^{(r)}(m;\mu,\alpha) = \mathbb{1}\left\{\widehat{\mu}_t - \mu \geqslant \widehat{\sigma}_t\cdot c_\alpha^{(o)}(\rho)\cdot b_t(m;\rho)\right\},\quad \varphi_t^{(\ell)}(m;\mu,\alpha) = \mathbb{1}\left\{\widehat{\mu}_t - \mu \leqslant -\,\widehat{\sigma}_t\cdot c_\alpha^{(o)}(\rho)\cdot b_t(m;\rho)\right\}, \tag{3}$$

where the critical values $c_\alpha^{(o)}(\rho)$ are chosen as the $(1-\alpha)$-quantile of $\sup_{y>0}(\rho(y)\cdot W(y))$ for a standard Wiener process $(W(y))_{y\in(0,\infty)}$. Unlike for a posteriori two-sided versus one-sided location tests, here, $c_\alpha(\rho) < c_{\alpha/2}^{(o)}(\rho)$: While realisations of random variables can only be either above the upper or below the lower quantile, a sample path of a process can first cross the lower and then still cross the upper boundary curve later.

The hypotheses in (2) are denoted as ordered (Lei and Fithian, 2016) or hierarchical (Rom and Holland, 1995).

**Corollary 3.3.** *Under the assumptions of Theorem 3.1 the FWER of simultaneously using the tests in (3) for* $H_0^{(r)}(\mu_1),\ \ldots, H_0^{(r)}(\mu_k)$ *with any* $\mu_1 < \cdots < \mu_k$, $k\in\mathbb{N}$, *is controlled asymptotically, i.e. for any* $\mu_X$ *it holds*

$$\limsup_{m\to\infty}\mathbb{P}_{\mu_X}\left(\varphi_t^{(r)}(m;\mu,\alpha) = 1 \text{ for some } \mu_j \geqslant \mu_X, j = 1,\ldots,k, \text{ and some } t\in\mathbb{N}\right) \leqslant \alpha.$$

*An analogous assertion holds for* $H_0^{(\ell)}(\mu_1),\ldots,H_0^{(\ell)}(\mu_k)$*, and even for simultaneously testing both left- and right-sided hypotheses, if in the tests in (3) the critical values* $c_\alpha^{(o)}(\rho)$ *are replaced by* $c_\alpha(\rho)$.

**Proof.** By construction of the tests, the probability in the corollary is dominated by $\mathbb{P}_{\mu_X=\mu}(\varphi_t^{(r)}(m;\mu_X,\alpha) = 1 \text{for some } t\in\mathbb{N})$, so that the result follows analogously to Theorem 3.1 from Remark 4.4. The left-sided result follows by symmetry. When simultaneously testing both left- and right-sided alternatives with the two-sided critical values, the probability of any false rejection at any time is upper-bounded by $\mathbb{P}_{\mu_X=\mu}(\varphi_t(m;\mu_X,\alpha) = 1 \text{ for some } t\in\mathbb{N})$ by construction. $\square$

Analogously to Theorem 3.2, the stopping guarantee under alternatives holds for these hierarchical one-sided tests.

## 4. Underlying asymptotic theory

The results of the previous two sections are based on the following limit theorem.

**Theorem 4.1.** *Let Assumptions $\mathcal{M}$ and $\mathcal{BC}$ be fulfilled, where instead of assuming that $\widehat{\sigma}_t > 0$ a.s. in $\mathcal{M}$, we may replace $1/\widehat{\sigma}_t$ below by $1/\widehat{\sigma}_t \mathbb{1}\left\{\widehat{\sigma}_t > 0\right\}$ with the convention that $0/0 = 0$. Then, it holds for any $\mathfrak{l}_m/m \to 0$ that*

$$\sup_{t\geqslant \mathfrak{l}_m}\left|\frac{\rho(t/m)}{\sqrt{m}}\cdot\frac{1}{\widehat{\sigma}_t}\cdot\sum_{j=1}^{t}\left(X_j-\mu_X\right)\right| \xrightarrow{\mathcal{D}} \sup_{y>0}|\rho(y)\cdot W(y)| = Z_\rho,$$

*where $(W(y))_{y\in(0,\infty)}$ denotes a standard Wiener process.*

In the previous Theorems 2.2 and 3.1 we have used $\mathfrak{l}_m = 1$. However, in practice, it can be beneficial to use e.g. $\mathfrak{l}_m = \log(m)$ or $\mathfrak{l}_m = \sqrt{m}$ to guarantee a somewhat more stable behaviour of $\widehat{\sigma}_t$ and avoid false positives for very early time points. This is different from conducting tests only after the burn-in period of $m$ observations, as waiting for $l_m$ observations does not change the above limit while waiting for $m$ observations, i.e. choosing $\rho(s) = 0$ for $s < 1$, does.

**Remark 4.2.** In the proof, we only need the observation sequence to fulfil a functional central limit theorem with asymptotic variance $\sigma^2$ as well as the two generalised Hájek–Rényi inequalities in (5) and (7) to hold. All three assertions follow from strong invariance principles as required in Waudby-Smith et al. (2024, Theorem 2.8). Such strong invariance principles have been proven for a variety of time series under weak nonparametric assumptions, where typically the asymptotic variance is given by the so-called long-run variance taking the autocovariances of all lags into account; see Aue and Kirch (2024, Section 3.2) for more details.

**Proof of Theorem 4.1.** By the assumptions on the sequence of variance estimators, possibly replacing $1/\widehat{\sigma}_t$ by $1/\widehat{\sigma}_t \cdot \mathbb{1}\left\{\widehat{\sigma}_t > 0\right\}$, it holds

$$(V_1)\quad \sup_{t>T_m}\left|\frac{1}{\widehat{\sigma}_t}-\frac{1}{\sigma_X}\right| \to 0 \quad a.s. \quad (\text{as } T_m\to\infty), \qquad\qquad (V_2)\quad \sup_{t\geqslant 1}\frac{1}{\widehat{\sigma}_t} = O(1) \quad a.s.$$

Without loss of generality, let $\mu_X = 0$. For any fixed $(v, V] \subset (0, e_\rho]$ it follows from $(V_1)$, the functional central limit theorem, the uniform continuity of $\rho(\cdot)$ on the domain $[v, V]$ and $\mathfrak{l}_m/m \to 0$ as $m \to \infty$

$$\sup_{mv<t\leqslant mV}\left|\frac{\rho\left(\frac{t}{m}\right)\mathbb{1}_{\{t>\mathfrak{l}_m\}}}{\sqrt{m}}\frac{1}{\widehat{\sigma}_t}\sum_{i=1}^{\min(t,mV)}X_i\right| = \sup_{v<y\leqslant V}\left|\frac{\rho(y)}{\sqrt{m}}\frac{1}{\sigma_X}\sum_{i=1}^{\min(\lfloor my\rfloor,mV)}X_i\right|\cdot(1+o_P(1)) \xrightarrow{\mathcal{D}} \sup_{v<y\leqslant V}|\rho(y)\,W(\min(y,V))|\,. \tag{4}$$

Let $\widetilde{\gamma}_1 \in (\gamma_1, 1/2)$. Then, the generalised Hájek–Rényi inequality

$$\sup_{1\leqslant t<m}\frac{1}{m^{\frac{1}{2}-\widetilde{\gamma}_1}t^{\widetilde{\gamma}_1}}\left|\sum_{i=1}^{t}X_i\right| = O_P(1) \text{ uniformly in } m \tag{5}$$

holds as by the classical Hájek–Rényi inequality (Hájek and Rényi, 1955) for any $C > 0$ uniformly in $m$

$$\mathbb{P}\left(\sup_{1\leqslant t<m}\frac{1}{m^{\frac{1}{2}-\widetilde{\gamma}_1}t^{\widetilde{\gamma}_1}}\left|\sum_{i=1}^{t}X_i\right| > C\right) \leqslant \frac{1}{C^2}\,\sigma_X^2\,\frac{1}{m^{1-2\widetilde{\gamma}_1}}\sum_{t=1}^{m}\frac{1}{t^{2\widetilde{\gamma}_1}} \leqslant \frac{1}{C^2}\,\frac{\sigma_X^2}{1-2\widetilde{\gamma}_1} = O\left(\frac{1}{C^2}\right).$$

Combining (5) with $(A_1)$ from Assumption $\mathcal{BC}$ and $(V_2)$, it follows for $v \to 0$ uniformly in $m$

$$\sup_{1\leqslant t\leqslant mv}\left|\frac{\rho\left(\frac{t}{m}\right)\mathbb{1}_{\{t>\mathfrak{l}_m\}}}{\sqrt{m}}\frac{1}{\widehat{\sigma}_t}\sum_{i=1}^{\min(t,mV)}X_i\right| \leqslant \sup_{0<s<v}s^{\widetilde{\gamma}_1}\rho(s)\cdot\sup_{1\leqslant t<m}\frac{1}{m^{\frac{1}{2}-\widetilde{\gamma}_1}t^{\widetilde{\gamma}_1}}\left|\sum_{i=1}^{t}X_i\right|\cdot\sup_{t\geqslant 1}\frac{1}{\widehat{\sigma}_t}\xrightarrow{P} 0. \tag{6}$$

Similarly, let $\widetilde{\gamma}_2 \in (\gamma_2, 1/2)$. Then, the generalised Hájek–Rényi inequality

$$\sup_{t>mV}\frac{(mV)^{\frac{1}{2}-\widetilde{\gamma}_2}}{t^{1-\widetilde{\gamma}_2}}\left|\sum_{i=\lfloor mV\rfloor+1}^{t}X_i\right| = O_P(1) \text{ uniformly in } m \text{ for any } V>0 \tag{7}$$

holds: Indeed, using similar proof techniques as in Frank (1966) and continuity from below, the classical Hájek–Rényi inequality can be extended to an unbounded domain, i.e.

$$\mathbb{P}\left(\sup_{t>mV}\frac{(mV)^{\frac{1}{2}-\widetilde{\gamma}_2}}{t^{1-\widetilde{\gamma}_2}}\left|\sum_{i=\lfloor mV\rfloor+1}^{t}X_i\right| > C\right) \leqslant \frac{\sigma_X^2}{C^2}\cdot\left(1+(mV)^{1-2\widetilde{\gamma}_2}\sum_{i=mV+1}^{\infty}\frac{1}{i^{2-2\widetilde{\gamma}_2}}\right) \leqslant \frac{\sigma_X^2}{C^2}\cdot\left(1+\frac{1}{1-2\widetilde{\gamma}_2}\right) = O\left(\frac{1}{C^2}\right).$$

Combining (7) with $(A_2)$ from Assumption $\mathcal{BC}$ and $(V_2)$ we get for $V \to \infty$ uniformly in $m$

$$\sup_{t\geqslant 1}\left|\frac{\rho\left(\frac{t}{m}\right)\mathbb{1}_{\{t>\mathfrak{l}_m\}}}{\sqrt{m}}\frac{1}{\widehat{\sigma}_t}\left(\sum_{i=1}^{\min(t,mV)}X_i-\sum_{i=1}^{t}X_i\right)\right| \leqslant \sup_{t\geqslant 1}\frac{1}{\widehat{\sigma}_t}\cdot\sup_{s>V}s^{1-\widetilde{\gamma}_2}\rho(s)\cdot\sup_{t>mV}\frac{(mV)^{\frac{1}{2}-\widetilde{\gamma}_2}}{t^{1-\widetilde{\gamma}_2}}\left|\sum_{i=\lfloor mV\rfloor+1}^{t}X_i\right|\cdot V^{\widetilde{\gamma}_2-\frac{1}{2}}\xrightarrow{P} 0. \tag{8}$$

By the law of the iterated logarithm for the Wiener process (Csörgő and Révész, 1981, Theorem 1.3.1) and $(A_2)$ from Assumption $\mathcal{BC}$, it holds for any $\upsilon > 0$

$$\sup_{s>\upsilon} |\rho(s)[W(s) - W(\min(s, V))]| \leqslant \sup_{s>V} |\rho(s)\, W(V)| + \sup_{s>V} |\rho(s)\, W(s)| \overset{P}{\longrightarrow} 0 \text{ for } V \to \infty. \tag{9}$$

Moreover, the law of the iterated logarithm (Csörgő and Révész, 1981, Theorem 1.3.1) in combination with Csörgő and Révész (1981, Lemma 1.3.3) and $(A_1)$ from Assumption $\mathcal{BC}$ yields for $\upsilon \to 0$

$$\sup_{0<s\leqslant\upsilon} |\rho(s)\, W(s)| = \sup_{t\geqslant 1/\upsilon} \left|\rho\left(\frac{1}{t}\right) W\left(\frac{1}{t}\right)\right| \overset{D}{=} \sup_{t\geqslant 1/\upsilon} \left|\frac{\rho\left(\frac{1}{t}\right) W(t)}{t}\right| \leqslant \sup_{t\geqslant 1/\upsilon} \left|\frac{W(t)}{t^{1-\widetilde{\gamma}_1}}\right| \cdot \sup_{t\geqslant 1/\upsilon} \left(\frac{1}{t}\right)^{\widetilde{\gamma}_1} \rho\left(\frac{1}{t}\right) \overset{P}{\longrightarrow} 0. \tag{10}$$

Stöhr (2019, Lemma B.2.) in combination with (4), (6) and (8)–(10) completes the proof. □

**Proposition 4.3.** *By choosing $\rho(s) = (1+s)^{\gamma_1+\gamma_2-1}/s^{\gamma_1}$ for some $0 \leqslant \gamma_1, \gamma_2 < 1/2$, it holds*

$$\sup_{y>0} |\rho(y) \cdot W(y)| \overset{D}{=} \sup_{0\leqslant x\leqslant 1} \frac{|B(x)|}{x^{\gamma_1}(1-x)^{\gamma_2}}$$

*for the limit distribution in Theorem 2.2, where $(B(y))_{y\in[0,1]}$ denotes a standard Brownian bridge.*

**Proof.** By substituting $x = \frac{y}{1-y}$ it holds

$$\sup_{y>0} |\rho(y)\, W(y)| = \sup_{x\in(0,1)} \left| \frac{\rho\left(\frac{x}{1-x}\right)}{1-x} (1-x)\, W\left(\frac{x}{1-x}\right)\right| \overset{D}{=} \sup_{x\in(0,1)} \frac{|B(x)|}{x^{\gamma_1}(1-x)^{\gamma_2}},$$

where the distributional equality follows from Csörgő and Révész (1981, Equation (1.4.5)). □

**Remark 4.4.** The statements of Theorem 4.1 and Proposition 4.3 remain true if all absolute values of quantities are replaced by the original quantities without absolute values. The proofs are analogous.

## 5. Discussion and outlook

In this paper, we propose a new class of boundary functions that can be used to obtain sharp uniform confidence sequences for the location model. These confidence sequences are dual to a class of sequential tests that can be used for safe testing and by construction even permit simultaneous testing of hierarchical one-sided hypotheses while still controlling the FWER among these multiple tests. The proposed nonparametric methodology is valid for a large class of random variables including time-dependent observations without the need to make specific assumptions on the underlying distributional family.

The results are based on asymptotic methodology that is strongly related to online monitoring methods for change points that have been first proposed by Chu et al. (1996), Horváth et al. (2004) and later extended in various directions, see Aue and Kirch (2024) for a recent survey on the subject. Several of these extensions are also particularly promising in the context of this paper, including extensions to the multivariate case (e.g. included in the setting of Kirch and Tadjuidje Kamgaing (2015)), more robust methodology for the location model such as sequential methods based on $M$-estimators (e.g. Koubková (2006), Chochola (2013), Kirch and Tadjuidje Kamgaing (2015)) or linear regression under various assumptions on the error structure (e.g. Chu et al. (1996), Horváth et al. (2004), Hušková and Koubková (2005), Aue et al. (2006)).

Furthermore, similar mathematical methods as in Hlávka et al. (2012) could be used to construct corresponding sequential goodness-of-fit tests. Extensions to two-sample testing where observations from both samples are regularly being taken are less straightforward even though the original sequential change point tests are extensions of two-sample tests but with the difference that the comparison is made between older and more recent data. Nevertheless, similar ideas could be used for mean-based methodology as well as to obtain sequential generalisations based on U-statistics including Wilcoxon-type methodology as has been used for change point monitoring by Kirch and Stoehr (2022).

## Funding

This work was financially supported by the Deutsche Forschungsgemeinschaft (DFG, German Research Foundation) - 314838170, GRK 2297 MathCoRe.

## Data availability

No data was used for the research described in the article.